\documentclass{article}
\usepackage{epsfig}
\usepackage{psfrag}

\usepackage{amsmath,amssymb,amscd,amsthm}
\title{
Continued fractions for some transcendental numbers}
\author{Andrew N.~W. Hone\thanks{{\bf Acknowledgments:} 
This work is supported by Fellowship EP/M004333/1 
from the Engineering and Physical Sciences Research Council.
The original inspiration came from Paul Hanna's  observations concerning the nonlinear recurrence sequences 
described in  \cite{curious}, which were communicated via the Seqfan mailing list. 
The author is grateful to Jeffrey Shallit for helpful correspondence on related matters.}\\
School of Mathematics,
Statistics and Actuarial Science,\\ University of
Kent,\\ Canterbury CT2 7NF, U.K.\\
\texttt{A.N.W.Hone@kent.ac.uk}
}

\newcommand{\br}{\begin{rem}}
\newcommand{\er}{\end{rem}}
\newcommand{\bex}{\begin{exa}}
\newcommand{\eex}{\end{exa}}
\newcommand{\bd}{\begin{Def}}
\newcommand{\ed}{\end{Def}}
\newcommand{\bt}{\begin{theorem}}
\newcommand{\et}{\end{theorem}}
\newcommand{\bl}{\begin{lemma}}
\newcommand{\el}{\end{lemma}}
\newcommand{\be}{\begin{equation}}
\newcommand{\ee}{\end{equation}}
\newcommand{\bea}{\begin{eqnarray}}
\newcommand{\eea}{\end{eqnarray}}

\newcommand{\adots}{\mathinner{\mkern2mu\raise1pt\hbox{.}\mkern2mu
\raise4pt\hbox{.}\mkern2mu\raise7pt\hbox{.}\mkern1mu}}

\newcommand{\beq}{\begin{equation}}
\newcommand{\eeq}{\end{equation}}
\newcommand{\bear}{\begin{array}}
\newcommand{\eear}{\end{array}}

\newcommand\ka{{\kappa}}
\newcommand\al{{\alpha}}

\newcommand\bM{\mathbf{M}}

\newtheorem{thm}{Theorem}[section]

\newtheorem{rem}[thm]{Remark}

\newtheorem{exa}[thm]{Example}

\newtheorem{lem}[thm]{Lemma}

\newenvironment{prf}{\trivlist \item [\hskip
\labelsep {\bf Proof:}]\ignorespaces}{\qed \endtrivlist}

\theoremstyle{remark}

\newcommand{\N}{{\mathbb N}}

\newcommand{\Z}{{\mathbb Z}}

\begin{document}

\maketitle

\begin{abstract}
We consider series of the form 
$$ 
\frac{p}{q} +\sum_{j=2}^\infty \frac{1}{x_j},  
$$  
where $x_1=q$ and the integer sequence $(x_n)$ satisfies a  
certain non-autonomous recurrence of second order, which 
entails that $x_n|x_{n+1}$ for $n\geq 1$. It is shown   
that the terms of the sequence, and multiples of the ratios of successive terms, 
appear interlaced in the continued fraction expansion of the sum of the 
series, which is a transcendental number. 

\vspace{.05in}

\vspace{.05in}
\end{abstract}

\section{Introduction}

\setcounter{equation}{0} 

In recent work \cite{curious}, we considered the integer sequence 
\beq\label{exs} 
1,1,2,12,936,68408496,342022190843338960032,\ldots
\eeq 
(sequence A112373 in Sloane's Online Encyclopedia of Integer Sequences),
which is generated from the initial values $x_0=x_1=1$ by the nonlinear recurrence relation 
\beq\label{orig} 
 x_{n+2} \, x_n = x_{n+1}^2(x_{n+1}+1), 
\eeq 
and proved some observations of Hanna, namely that the sum 
\beq\label{sumr} 
\sum_{j=1}^\infty \frac{1}{x_j} 
\eeq 
has the continued fraction expansion 
\beq\label{pat}
[x_0; y_0,x_1,y_1,x_2,\ldots, y_{j-1},x_{j},\ldots],  
\eeq
where $y_j=x_{j+1}/x_j\in\N$ and we use the notation 
$$ 
[a_0; a_1,a_2,a_3,\ldots,a_n,\ldots ]=a_0+\cfrac{1}{a_1+\cfrac{1}{a_2+\cfrac{1}{a_3+\ldots\cfrac{1}{a_n+\ldots}}}} 
$$ 
for continued fractions. 
Furthermore, we generalized this result by obtaining the explicit 
continued fraction expansion for the sum of reciprocals (\ref{sumr}) 
in the case of a sequence $(x_n)$ 
generated by a nonlinear recurrence of the form 
\beq\label{form} 
x_{n+1} \, x_{n-1} = x_{n}^2\, F(x_{n}), 
\eeq 
with $F(x)\in\Z_{\geq 0}[x]$ and $F(0)=1$; so (\ref{orig}) 
corresponds to the particular case $F(x)=x+1$. 

All of the recurrences (\ref{form}) exhibit 
the Laurent phenomenon \cite{fz}, and starting from $x_0=x_1=1$ they generate a sequence of positive integers satisfying $x_n|x_{n+1}$. The latter fact means that the sum (\ref{sumr}) is an Engel series 
(see Theorem 2.3 in Duverney's book \cite{duve}, for instance). 

The purpose of this note is to present a further generalization of the 
results in \cite{curious}, by considering a sum 
\beq \label{gsum} 
S=
\frac{p}{q}+\sum_{j=2}^\infty \frac{1}{x_j},  
\eeq 
with the terms $x_n$ satisfying the recurrence 
\beq\label{nform} 
x_{n+1} \, x_{n-1} = x_{n}^2\, (z_nx_{n}+1), 
\eeq 
for $n\geq 2$, where $(z_n)$ is a sequence of positive integers,
$x_1=q$, and $x_2$ is specified suitably. Observe that, in contrast to 
(\ref{form}), 
the recurrence (\ref{nform}) can be viewed as a non-autonomous
dynamical system for $x_n$, because the coefficient $z_n$ can vary independently (unless 
it is taken to be $G(x_n)$,  for some function $G$).  The same argument as used in \cite{curious}, based on Roth's theorem, shows the transcendence of any  number $S$ 
defined by a sum of the form (\ref{gsum}) with such a sequence $(x_n)$.  

\section{The Main Result}

\setcounter{equation}{0} 

We start with a rational number written in lowest terms as $p/q$, and suppose that the continued fraction of this number is given as 
\beq\label{pqcf} 
\frac{p}{q}= [a_0; a_1,a_2,a_3,\ldots,a_{2k}]
\eeq  
for some $k\geq 0$. Note that, 
in accordance with a comment on p.230 of \cite{sh3}, there is no loss of generality in assuming 
that the index of the final coefficient is even. For the convergents we 
denote numerators and denominators by $p_n$ and $q_n$, respectively, and 
use the correspondence between matrix products and continued fractions, 
which says that  
\beq\label{matid} 
\bM_n:= 
\left(\begin{array}{cc} p_n & p_{n-1} \\ q_n & q_{n-1} \end{array}\right)=
\left(\begin{array}{cc} a_0 & 1 \\ 1 & 0 \end{array}\right)
\left(\begin{array}{cc} a_1 & 1 \\ 1 & 0 \end{array}\right) 
\ldots 
\left(\begin{array}{cc} a_n & 1 \\ 1 & 0 \end{array}\right),
\eeq 
yielding the determinantal identity 
\beq\label{detid}\det \bM_n =  
p_nq_{n-1}-p_{n-1}q_n = (-1)^{n+1}.
\eeq 
Now for a given sequence $(z_n)$ of positive integers, we define a new sequence  $(x_n)$ by 
\beq \label{xdef} 
x_1=q, \qquad x_{n+1}=x_ny_{n-1}(x_nz_n+1) \quad \mathrm{for} \quad n\geq 1, 
\eeq 
where 
\beq\label{ydef} 
y_0 = q_{2k-1}+1, \qquad y_n=\frac{x_{n+1}}{x_n}\quad \mathrm{for} \quad n\geq 1. 
\eeq 
It is clear from (\ref{xdef}) and (\ref{ydef}) that $(x_n)$ is an increasing sequence of positive integers such that $x_n|x_{n+1}$ for all $n\geq 1$; $(y_n)$ also consists of positive integers, and is an increasing  sequence as well.  
The recurrence  (\ref{nform}) for $n\geq 2$ follows immediately from 
(\ref{xdef}) and (\ref{ydef}).
\begin{thm}\label{main}  
The partial sums of (\ref{gsum}) are given by   
$$ 
S_n:=
\frac{p}{q}+\sum_{j=2}^n \frac{1}{x_j} 
=[a_0;a_1,\ldots,a_{2(k+n-1)}]   
$$
for all $n\geq 1$, where the coefficients appearing after $a_{2k}$ are 
$$ 
a_{2k+2j-1} = y_{j-1}z_{j}, \qquad a_{2k+2j}=x_j \quad 
for \quad j\geq 1.
$$ 
\end{thm} 
\begin{prf}
For $n=1$, $S_1$ is just (\ref{pqcf}), and we note that 
$q_{2k-1}=y_0-1$ and $q_{2k}=q=x_1$. 
Proceeding by induction, we suppose that 
$q_{2k+2n-3}=y_{n-1}-1$ and $q_{2k+2n-2}=x_n$, and calculate the 
product 
$$ 
\bear{l}
\bM_{2k+2n}=\bM_{2k+2n-2} \left(\begin{array}{cc} a_{2k+2n-1} & 1 \\ 1 & 0 \end{array}\right)
\left(\begin{array}{cc} a_{2k+2n} & 1 \\ 1 & 0 \end{array}\right) \\ 
\qquad\qquad \, = 
\bM_{2k+2n-2}
\left(\begin{array}{cc} y_{n-1}z_{n} & 1 \\ 1 & 0 \end{array}\right)
\left(\begin{array}{cc} x_n & 1 \\ 1 & 0 \end{array}\right) \\ 
\qquad \qquad\, = \left(\begin{array}{cc} p_{2k+2n-2}& p_{2k+2n-3} \\ q_{2k+2n-2} & q_{2k+2n-3} \end{array}\right) 
\left(\begin{array}{cc} x_ny_{n-1}z_n &  y_{n-1}z_n \\ x_n & 1 \end{array}\right).
\eear 
$$ 
By making use of (\ref{xdef}) and (\ref{ydef}),
this gives $p_{2k+2n}=(x_ny_{n-1}z_n+1)p_{2k+2n-2}+x_n p_{2k+2n-3}$,  
$$\bear{l}
q_{2k+2n-1} = y_{n-1}z_n\, q_{2k+2n-2}+ q_{2k+2n-3}
= x_n y_{n-1}z_n +y_{n-1}-1 \\
\qquad \qquad\,
=\frac{x_{n+1}}{x_n} -1 
= y_n-1,
\eear  
$$  
and 
$$ \bear{l}
q_{2k+2n}=(x_ny_{n-1}z_n+1)q_{2k+2n-2}+x_n q_{2k+2n-3}\\ 
\qquad \quad\,
=  (x_ny_{n-1}z_n+1)x_n + x_n (y_{n-1}-1)=x_{n+1},
\eear 
$$
which are the required denominators for the $(2k+2n-1)$th and 
$(2k+2n)$th convergents. Thus we have 
$$
S_{n+1}=S_n + \frac{1}{x_{n+1}} =\frac{p_{2k+2n-2}}{q_{2k+2n-2}} 
+\frac{1}{q_{2k+2n}} 
= \frac{1}{q_{2k+2n}}\left(\frac{x_{n+1}}{x_n}p_{2k+2n-2}+1\right). 
$$
From (\ref{detid}) and (\ref{xdef}), the bracketed expression above can be rewritten 
as 
$$\bear{l} 
\Big(y_{n-1}(x_nz_n+1) -q_{2n+2k-3}\Big)p_{2k+2n-2} +q_{2n+2k-2}p_{2k+2n-3} 
\\ 
\qquad = \Big(y_{n-1}(x_nz_n+1) -y_{n-1}+1\Big)p_{2k+2n-2} +x_np_{2k+2n-3},  
\eear 
$$ 
giving 
$$ 
S_{n+1} = \frac{1}{q_{2k+2n}}\Big( (x_ny_{n-1}z_n+1)p_{2k+2n-2}+x_n p_{2k+2n-3}\Big)= \frac{p_{2k+2n}}{q_{2k+2n}},
$$
which is the required result.
\end{prf} 

Upon taking the limit $n\to\infty$ we obtain the infinite continued fraction 
expansion for the sum $S$, which is clearly irrational. To show that $S$ is transcendental, we need the following growth estimate for $x_n$: 
\begin{lem}\label{lemma} The terms of a sequence defined by (\ref{xdef}) 
satisfy 
$$ x_{n+1}>x_n^{5/2}
$$  
for all $n\geq 3$.
\end{lem} 
\begin{prf} Since $(x_n)$ is an increasing sequence, the recurrence relation (\ref{nform}) gives 
$$ 
x_{n+1}>\frac{x_n^3}{x_{n-1}}>x_n^2 
$$ 
for $n\geq 2$. Hence $x_{n-1}<x_n^{1/2}$ for $n\geq 3$, and putting this back into the first inequality above yields 
$x_{n+1}>x_n^3/x_n^{1/2}=x_n^{5/2}$, as required.  
\end{prf} 

The preceding growth estimate for $x_n$ means that $S$ can be well approximated by rational numbers. 
\begin{thm}\label{trans} The sum 
$$ 
S= \frac{p}{q}+\sum_{j=2}^\infty \frac{1}{x_j}
=[a_0;a_1,\ldots,a_{2k},y_0z_1,x_1,y_1z_2,\ldots, y_{j-1}z_j,x_j, \ldots]
$$ 
is a transcendental number.
\end{thm}  
\begin{prf} This is the same as the proof of Theorem 4 in \cite{curious},  
which we briefly outline  here.  Let $P_n=p_{2k+2n-2}$ and $Q_n=q_{2k+2n-2}$.
Approximating the irrational number $S$ by the partial sum $S_n=P_n/Q_n$, then using Lemma \ref{lemma} and a comparison with a geometric sum, gives the upper bound 
$$
\left|S-\frac{P_n}{Q_n} \right|=\sum_{j=n+1}^\infty \frac{1}{x_j} <
\frac{1}{x_{n}^{5/2-\epsilon}}=\frac{1}{Q_n^{5/2-\epsilon}}
$$
for any $\epsilon>0$, whenever $n$ is sufficiently large. 
Roth's theorem \cite{roth} (see also chapter VI in \cite{cassels}) says that, for an arbitrary fixed $\ka>2$, an irrational algebraic number $\al$ has  only finitely many rational approximations $P/Q$ for which 
$
\left|\al -\frac{P}{Q}\right|<\frac{1}{Q^{\ka}};
$
so $S$ is transcendental. 
\end{prf}

For other examples of transcendental numbers whose continued fraction 
expansion is explicitly known, see \cite{ds} and references therein. 

\section{Examples}

\setcounter{equation}{0} 

The autonomous recurrences (\ref{form}) considered in \cite{curious}, where the polynomial $F$ has positive integer coefficients and $F(0)=1$,  give 
an infinite family of examples. In that case, one has $p=1$ and 
$x_1=q=1$, so that $k=0$, $y_0=1$ and $z_n=(F(x_n)-1)/x_n$. More generally, one could take $z_n=G(x_n)$ for any non-vanishing 
arithmetical function $G$. 

In general, it is sufficient to take the initial term in (\ref{gsum}) lying in the range $0<p/q\leq 1$, since going outside this range only alters the value of $a_0$. 
As a particular example\footnote{The published version of this example in Monatsh. Math. contains errors.}, we take 
$$ 
\frac{p}{q}=\frac{6}{7}=[0;1,6], \qquad z_n=n   \quad \mathrm{for}  \quad n\geq 1,
$$  
so that $k=1$, and $q_{1}=1$ which gives $y_0=2$. Hence $x_1=7$, 
$x_2=112$, and the sequence $(x_n)$ continues  with  
$$ 
403200,1755760043520000,53695136666462381094317154204367872000000,\ldots .
$$ 
The sum $S$ is the transcendental number 
$$ 
\frac{6}{7}+\frac{1}{112}+\frac{1}{403200}+ \frac{1}{1755760043520000}+\ldots \approx 
0.86607390873015929971, 
$$ 
with continued fraction 
expansion 
$$ 
[0;1,6,2,7,32,112,10800,403200,17418254400,1755760043520000,\ldots]. 
$$


\end{document}